\documentclass[12pt]{article}
\usepackage[T1]{fontenc}
\usepackage[utf8]{inputenc}
\usepackage{hyperref}
\usepackage{amssymb}
\usepackage{amsmath}
\usepackage{array}
\usepackage{graphicx}
\usepackage[justification=centering]{caption}
\usepackage{chngcntr}
\counterwithin{figure}{section}

\title{Parallel-in-time method for calculation of long-range electrostatic interactions}
\author{Jana Paz\'urikov\'a \and Lud\v ek Matyska}
\date{\today}

\begin{document}

 \maketitle
 \section{Introduction}
  We are designing and developing the algorithms for molecular dynamics (MD) simulations with large systems ($\sim 10^{10}$ particles) over reasonably long time (tens of $\mu$s). The motivation came from a request from computational chemists at Faculty of Science, Masaryk University. Dr V\'acha has simulated the virus passing through the cell membrane, however only with a coarse-grained model where one bead represents tens of atoms (see \url{http://www.youtube.com/watch?v=Pkq-w2xhczw&feature=youtu.be}) \cite{Vacha2011}. Single atom precision and explicitly modelled water molecules would facilitate the deeper understanding of the process. The fine-grained model we want to simulate has following features:
  \begin{itemize}
    \item the size---the number of particles---$6-210.10^6$ atoms (with differently sized viruses and different non-atomic resolution on the water solution);
    \item the box---the dimensions---$30\text{x}65\text{x}65\text{ nm}$;
    \item the time---the number of algorithm's iterations---$50\ \mu s \sim25.10^9$ steps with $2\text{ fs}$ timestep.
   
  \end{itemize}
  
  Several scientific areas impact the research of large scale MD simulations:
  \begin{itemize}
  \item biology and chemistry pose the problem: what is happening at atomic level when the virus is passing though the cell membrane?
  \item computational chemistry offers the model to simulate the problem---molecular dynamics; empirical equations and parameters to approximate interactions between atoms; algorithms to overcome the main bottleneck---the electrostatic interactions;
  \item physics frames the world of simulation: in molecular dynamics the atoms behave according to Newton's second law of motion;
  \item mathematics formulates Newton's second law of motion as partial differential equations and provides numerical methods to solve them;
  \item we, as computer scientists, aim to develop and implement the highly parallelizable algorithm capable of simulating such large system for such a long time.
  \end{itemize}

  The report has five more sections. In the second section, we shortly introduce the molecular dynamics and the main bottleneck of the simulation---the electrostatic interactions. In the third section, we explain the algorithm for calculation of electrostatic interactions---the multilevel summation method. In the fourth section, we explain parallel-in-time numerical method---the parareal method. In the fifth section, we present our novel method, the combination of multilevel summation method and parareal method. In the last section, we propose its possible improvements and list open problems where more explicit mathematical analysis is needed.
 
 \section{Molecular Dynamics}
 Molecular dynamics observes and experiments with realistic model of molecu\-les in order to understand chemical or biological processes and predict macroscopic properties by detailed knowledge of particles's movements caused by all-to-all interactions \cite{Jensen2007}. The atoms interact if they are bonded---due to bond stretching, angle bending, dihedral torsions---but also if they are not bonded---due to van der Waals and Coulomb/electrostatic interactions. The potential (or energy) between the particles is expressed as the force that changes the particles's positions. The movement is governed by Newton's second law of motion commonly known as $\mathbf{F} = m\mathbf{a}$ where $\mathbf{F}$ is the force, $m$ is the mass of the atom, $\mathbf{a}$ is the acceleration.

A general simulation algorithm first takes input data---partial charges $q_i$, positions $\mathbf{r}_i$ and velocities $\mathbf{v}_i$ for all particles in the system---and then iteratively repeats the following steps:
\begin{enumerate}
 \item {calculate the potential and forces}
 \begin{equation}
  U_{all} = U_{bonded} + U_{non-bonded} = U_{bond} + U_{angle} + U_{dihed} + U_{Waals} + U_{Coulomb}
 \end{equation} 
 where $U_*$ is the potential due to $*$;
 
 \begin{equation}
  \mathbf{F}_i =-\frac{\partial U_{all}}{\partial \mathbf{r}_i}
 \end{equation}
where $\mathbf{F}_i$ is the force exerted on atom $i$;

 \item {move particles}
 
 \begin{equation}
 \label{eq:pde}
  \mathbf{F}_i =- \frac{\partial U_{all}}{\partial \mathbf{r}_i} = m_i \mathbf{a}_i = m_i \frac{\partial \mathbf{v}_i}{\partial t} = m_i\frac{\partial^2 \mathbf{r}_i}{\partial t^2};
 \end{equation}
 
 \item update time, optionally generate output.
\end{enumerate}

The output of the simulation includes the trajectories of particles, forces and energy of the system. Properties of the system such as average potential energy, temperature or the viscosity of the liquid are processed from output data by using further physical and chemical equations and applying statistical methods.

Each $U_*$ is approximated by empirical equations and parameters based on the atoms involved. The partial differential equation (\ref{eq:pde}) is solved by numerical methods, most common are leapfrog and velocity Verlet integration schemes.

The bottleneck of the computation, long-range electrostatic interactions $U_{Coulomb}$ are calculated by Coulomb's law between almost every pair of atoms (thus with intrinsic quadratic complexity) as
\begin{equation}
 \label{eq:coulomb}
 U_{Coulomb} = \sum_i \sum_j \frac{q_i q_j}{|\mathbf{r}_i - \mathbf{r}_j|}.
\end{equation}

The research in molecular dynamics focuses intensively on methods for Coulomb's law approximation \cite{Koehl2006}. The simplest method---simple cutoff method---calculates the interactions by Coulomb law but only for atoms that are within cutoff distance. The mwthod is inaccurate and introduces artificial physical artefacts are the edge of cutoff sphere. Wolf summation method deals with most of these disadvantages. 

All more sophisticated and more accurate methods interpolate the charges onto the grid and then they calculate the potential by
\begin{itemize}
 \item Fourier transform---Ewald sum \cite{ewald21}, Particle-mesh Ewald \cite{darden93}, Smooth Particle-mesh Ewald \cite{essmann95}, Gaussian split Ewald \cite{shan05};
 \item hierarchical division of the space---Barnes-Hut method \cite{barneshut86}, Fast multipole method \cite{greengard87};
 \item multigrid methods \cite{Stuben1982, Briggs2000}.
\end{itemize}

Most software packages for MD simulation implement the methods based on Fourier transform. 
However, they do not scale well in number of processors used due to many-to-many communication pattern. 

Methods based on hierarchical division of space stems from astronomical N-body simulations that resemble MD---instead of electrostatic potential they calculate gravitational interactions. Fast multipole method has $\mathcal{O}(N)$ complexity, however with a large multiplicative constant. Moreover, the implementation is rather difficult.

Multigrid methods, a mathematical approach to solve partial differential equations, apply local process to multiple scales (grids) of the problem. They obtain an initial approximation of results on the finest grid. As the error is smooth, the corrections can be calculated on coarser grids by recursive relaxation (e.g.~Jacobi relaxation method) and restriction to coarser grids. The correction then prolongates to finer and finer grids. Their large multiplicative constant for $\mathcal{O}(N)$ complexity and iterative nature make them inefficient for massively parallel implementation. 

Nevertheless, Hardy in his dissertation \cite{Hardy2006} took the multiscale approach and developed the multilevel summation method (MSM) that overcomes mentioned disadvantages of common methods. Therefore we decided to include it in our simulation. 

 \section{Multilevel Summation Method}
 \label{sec:multilevel}
   Hardy in his dissertation \cite{Hardy2006} developed the multilevel summation method (MSM) for calculation of long-range interactions with thorough mathematical background, performance assessments, accuracy analysis and implementation suggestions. The method divides the calculation of the potential onto multiple grids. To keep reasonable accuracy, more and more slowly varying parts of the potential are calculated on coarser and coarser grids.

  

  Multilevel summation method hierarchically interpolates smoothed potential onto multiple grids and then sums it up \cite{Hardy2009}. Electrostatic potential $U_{Coulomb}$ depends only on charges that are given as simulation's input and the reciprocal of the distance between atoms, i.e. $\frac{1}{|\mathbf{r}_{j}-\mathbf{r}_{i}|}$. We can rewrite the expression so that the sum gives the same result, but its parts correspond to more and more slowly varying parts of potential.
  \begin{equation}
  \label{eq:split}
  \begin{gathered}
   \frac{1}{|\mathbf{r}_{j} - \mathbf{r}_{i}|} = \left( \frac{1}{|\mathbf{r}_{j} - \mathbf{r}_{i}|} - g_a(\mathbf{r}_i, \mathbf{r}_j)\right) + g_a(\mathbf{r}_i, \mathbf{r}_j) \\
   g_a(\mathbf{r}_i, \mathbf{r}_j) = \frac{1}{a}\gamma \left( \frac{|\mathbf{r}_j-\mathbf{r}_i|}{a} \right) \\
   g_a(\mathbf{r}_i, \mathbf{r}_j) = \left( g_a(\mathbf{r}_i, \mathbf{r}_j) - g_{2a}(\mathbf{r}_i, \mathbf{r}_j) \right) + g_{2a}(\mathbf{r}_i, \mathbf{r}_j)
   \end{gathered}
  \end{equation}
where $a$ is the cutoff distance and $\gamma$ is the smoothing function.

We can rewrite these equations, instead of subscript $a$ representing the cutoff distance, we can use superscript $k$ representing the grid level. We will use the superscript notation from now on.
  
\begin{equation}
  \begin{gathered}
      \frac{1}{|\mathbf{r}_{j}-\mathbf{r}_{i} |} = (g^{*} + g^{0} + g^{1} + ... + g^{l-2} + g^{l-1})(\mathbf{r}_i, \mathbf{r}_j) \\
   g^*(\mathbf{r}_i, \mathbf{r}_j) = \frac{1}{|\mathbf{r}_{j}-\mathbf{r}_{i}|}-\frac{1}{a}\gamma(\frac{|\mathbf{r}_{j}-\mathbf{r}_{i}|}{a}) \\
   g^k(\mathbf{r}_i, \mathbf{r}_j) = \frac{1}{2^ka}\gamma(\frac{|\mathbf{r}_{j}-\mathbf{r}_{i}|}{2^ka})-\frac{1}{2^{k+1}a}\gamma(\frac{|\mathbf{r}_{j}-\mathbf{r}_{i}|}{2^{k+1}a}) \text{ for } k=0..l-2\\
   g^{l-1}(\mathbf{r}_i, \mathbf{r}_j) = \frac{1}{2^{l-1}a} \gamma(\frac{|\mathbf{r}_{j}-\mathbf{r}_{i}|}{2^{l-1}a}) 
  \end{gathered}
\end{equation}
where $l$ is number of grids (grid levels $0..l-1$).
 
 $\gamma$ is an unparameterized smoothing of the function $\frac{1}{\rho}$ chosen so that the first part of the first equation \ref{eq:split}, $\left( \frac{1}{|\mathbf{r}_{j} - \mathbf{r}_{i}|} - g_a(\mathbf{r}_i, \mathbf{r}_j)\right)$, vanishes after the cutoff distance ($a$ for the finest grid, $2^k a$ for grids of level $k$). The second part, $g_a(\mathbf{r}_i, \mathbf{r}_j)$, has partial derivatives and is slowly varying, which depends mainly on function $\gamma$. This function is usually Taylor expansion. 
 \begin{equation}
  \gamma(\rho) = \left\{ \begin{array}{lr} \frac{15}{8}-\frac{5}{4}\rho^2 + \frac{3}{8}\rho^4, & \rho < 1,\\ 1/\rho, & \rho \geq 1 \end{array}\right. 
 \end{equation}
 
 We have divided the potential into several, more and more slowly varying parts. How do we put them on multiple grids and approximate the solution? For each grid with spacing $2^kh$ and grid points $r_{\mu}^{k}$ \footnote{Letters $\mu$ and $\nu$ represent grids or grid points, $k$ is the grid level, $h$ is the spacing of the finest grid.}, we define nodal basis functions $\phi_{\mu}^k$ with local support (the function in non-zero only in the close surrounding of the grid point) as \begin{equation} \phi_{\mu}^k = \Phi\left(\frac{x-x_{\mu}^k}{2^kh}\right) \Phi\left(\frac{y-y_{\mu}^k}{2^kh}\right) \Phi\left(\frac{z-z_{\mu}^k}{2^kh}\right),\end{equation} where $\Phi(\epsilon)$ can be for example cubic interpolating polynomial. 
 
 Interpolation operator $\mathcal{I}^k$ interpolates each $g^k(\mathbf{r}_i, \mathbf{r}_j)$ to the grid by 
 \begin{equation}
 \begin{gathered}
  \mathcal{I}^k g(\mathbf{r}_i, \mathbf{r}_j) = \sum_{\mu} \sum_{\nu} \phi_{\mu}^k(\mathbf{r}_i) g(\mathbf{r}_{\mu}^k, \mathbf{r}_{\nu}^k) \phi_{\nu}^k(\mathbf{r}_j)\\
  k=0,1...l-1
  \end{gathered}
 \end{equation}

 Then, we approximate $\frac{1}{|\mathbf{r}_i-\mathbf{r}_j|}$ as
 \begin{equation}
  \begin{gathered}
    \frac{1}{|\mathbf{r}_i-\mathbf{r}_j|} = (g^* + g^0 + g^1 + ... + g^{l-2} + g^{l-1})(\mathbf{r}_i, \mathbf{r}_j) \\
    \approx (g^* + \mathcal{I}^0(g^0 + \mathcal{I}^1(g^1 + ... + \mathcal{I}^{l-2}(g^{l-2} + \mathcal{I}^{l-1}g^{l-1})...)))(\mathbf{r}_i, \mathbf{r}_j)
  \end{gathered}
 \end{equation}

 With approximated reciprocal distance, we can calculate the potential as
 \begin{equation}
  U_i \approx \frac{1}{4\pi \epsilon_0} (u_i^{short} + u_i^{long})
 \end{equation}

 The short-range part within cutoff is $u_i^{short} = \sum_{j} g^*(\mathbf{r}_{i}, \mathbf{r}_{j})q_j$. Long-range part is recursively divided between two parts---one within the cutoff and then calculated on the current grid; and the second representating even more slowly varying potential and then calculated on coarser grids. 
 The method for calculating the long-range part goes as follows \cite{Hardy2009}:
 \begin{itemize}
  \item \textbf{anterpolation}---puts point charges onto the grid $q_{\mu}^0 = \sum_{j}\phi_{\mu}^0(\mathbf{r}_j)q_{j}$;
  \item recursively for $k=0,1,...,l-2$
   \begin{itemize}
     \item \textbf{restriction}---approximates charges onto coarser grid \\$q_{\mu}^{k+1} = \sum_{\nu} \phi_{\mu}^{k+1}(\mathbf{r}_{\nu}^k)q_\nu^k$;
     \item \textbf{lattice cutoff}---calculates the part of potential corresponding to the grid $u_{\mu}^{k, cutoff} = \sum_{\nu}g^k(\mathbf{r}_{\mu}^k, \mathbf{r}_{\nu}^k) q_{\nu}^k$;
   \end{itemize}
  \item \textbf{top level}---calculates the most slowly varying part of the potential corresponding to the coarsest grid $u_{\mu}^{l-1} = \sum_{\nu} g^{l-1}(\mathbf{r}_\mu^{l-1}, \mathbf{r}_{\nu}^{l-1})q_{\nu}^{l-1}$;
  \item \textbf{prolongation}---recursively backwards for $k=l-2, ..., 1, 0$ adds up the parts of potential corresponding to the grids $u_\mu^k = u_\mu^{k, cutoff} + \sum_{\nu} \phi_{\nu}^{k+1}(\mathbf{r}_{\mu}^k)u_{\nu}^{k+1}$;
  \item \textbf{interpolation}---puts grid potential off grid $u_i^{long} = \sum_\mu \phi_\mu^0(\mathbf{r}_i)u_\mu^0$.
 \end{itemize} 
 
  Multilevel summation method exhibits several advantages over common methods for calculation of electrostatic interactions \cite[ch.1.2]{Hardy2006}. Apart from fast multipole method (FMM), it calculates continuous forces, smooth potential and therefore lower accuracy suffices for stable dynamics and conservation of energy. However, MSM does not conserve linear momentum. When compared to particle-mesh Ewald method (PME), the most obvious difference is only $\mathcal{O}(N)$ asymptotic complexity. Nevertheless, $\mathcal{O}(N.\log(N))$ is not much slower in practice and multiplicative constants delete performance results (actually, PME runs faster than FMM). Multilevel summation method scales better in parallel computation, uses multiple-time-step integration scheme more efficiently and communicates less than PME \cite[ch.7.4]{Hardy2006}. Apart from multigrid methods, MSM is not iterative therefore it scales better in parallel computation and does not require more than one 
global communication exchange. 
 
\begin{figure}[htb]
 \begin{center}
 \caption{Example of multilevel summation algorithm.}
 \label{fig:multilevel}
  \includegraphics[width=0.8\textwidth]{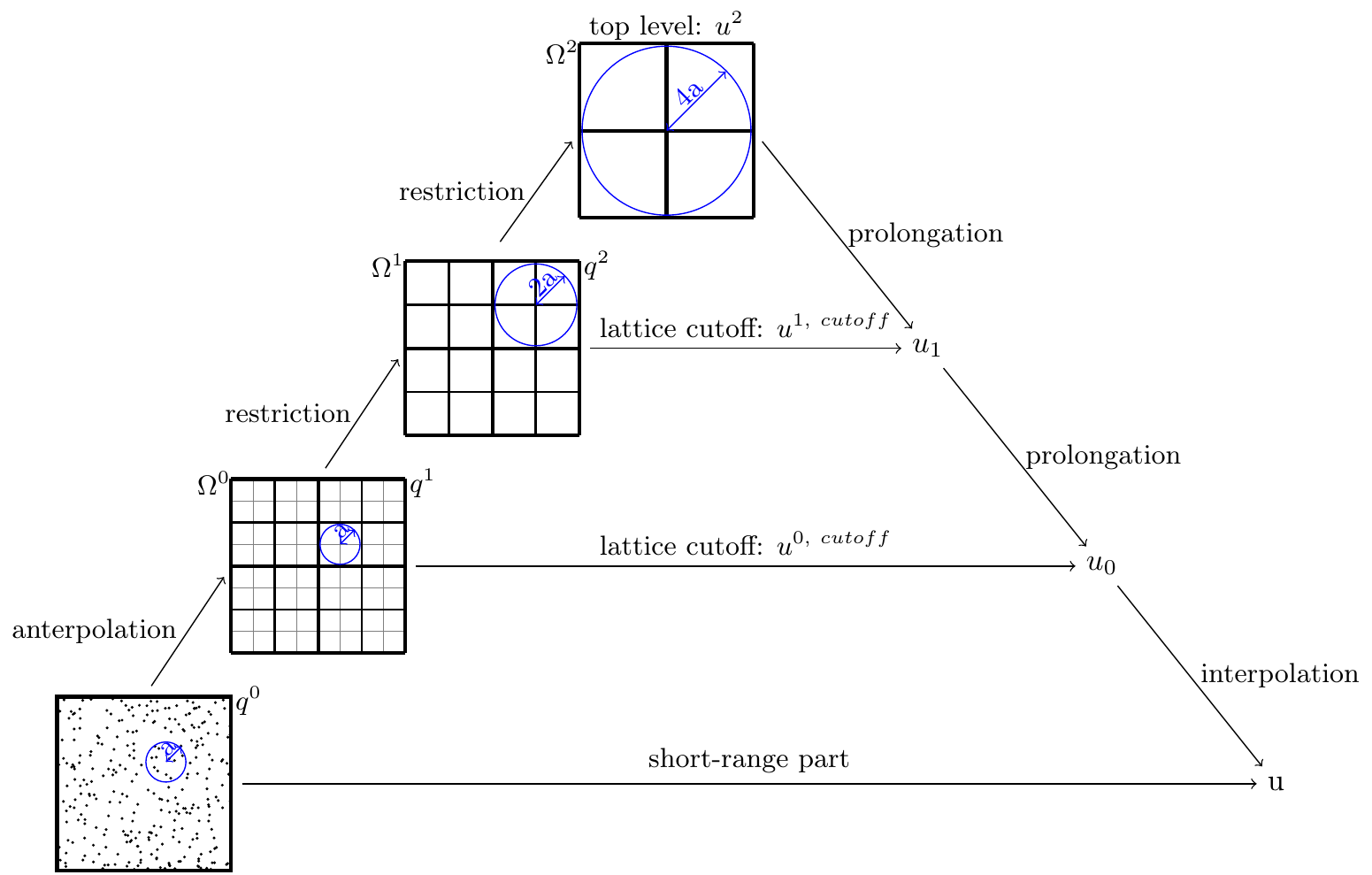}
  \end{center}
 \end{figure}

 The figure \ref{fig:multilevel} shows a specific example. We have three grids, $\Omega^0\text{, }\Omega^1\text{, }\Omega^2$, each with double spacing of the previous one. The method goes as:
 \begin{enumerate}
  \item We calculate short-range part of the potential $u_i^{short} = \sum_{j} g^*(\mathbf{r}_{i}, \mathbf{r}_{j})q_j$ with modified Coulomb's law (without $1/{4\pi\epsilon_0}$) for all atoms $j$ that are apart no more than cutoff distance $a$.
  \item We anterpolate the charges of all atoms to the finest grid. For each grid point, we take all atoms $j$ (with position $\mathbf{r}_j$ and charge $q_j$) that are close to the grid point (within local support of $\phi^0$) and calculate the grid point's charge as $q^0 = \sum_j \phi^0(\mathbf{r}_j)q_j$.
  \item We restrict charges from grid $\Omega^0$ to coarser grid $\Omega^1$ by $q^1=\sum_j \phi^1(\mathbf{r}^0_{j})q_j^0$. We calculate lattice cutoff from grid $\Omega^0$ -- the part of grid potential in grid point $\mathbf{r}^0$ as $u^{0,cutoff}=\sum_j g^0(\mathbf{r}^0, \mathbf{r}_j^0)q_j^0$ where $g^0(\mathbf{r}^0, \mathbf{r}_j^0) = \frac{1}{a} \gamma\left(\frac{|\mathbf{r}_j^0-\mathbf{r}^0|}{a}\right)-\frac{1}{2a} \gamma\left(\frac{|\mathbf{r}_j^1-\mathbf{r}^1|}{2a}\right)$. 
  \item We restrict charges from grid $\Omega^1$ to the coarsest grid $\Omega^2$ by $q^2=\sum_j \phi^2(\mathbf{r}^1_{j})q_j^1$. We calculate lattice cutoff from grid $\Omega^1$ -- the part of grid potential in grid point $\mathbf{r}^1$ as $u^{1,cutoff}=\sum_j g^1(\mathbf{r}^1, \mathbf{r}_j^1)q_j^1$ where $g^1(\mathbf{r}_1, \mathbf{r}_j^1) = \frac{1}{2a} \gamma\left(\frac{|\mathbf{r}_j^1-\mathbf{r}^1|}{2a}\right)-\frac{1}{4a} \gamma\left(\frac{|\mathbf{r}_j^2-\mathbf{r}^2|}{4a}\right)$. 
  \item At the top level grid $\Omega^2$, we calculate the most slowly varying part of the potential by $u^2 = \sum_j g^2(\mathbf{r}^2, \mathbf{r}^2_j)q_j^2$ where $g^2(\mathbf{r}^2, \mathbf{r}_j^2) = \frac{1}{4a} \gamma\left(\frac{|\mathbf{r}_j^2-\mathbf{r}^2|}{4a}\right)$.
  \item We prolongate the part of the potential from grid $\Omega^2$ to $\Omega^1$ by $u^1 = u^{1, cutoff} + \sum_j \phi^2(\mathbf{r}^1)u^2_j$.
  \item We prolongate the part of the potential from grid $\Omega^1$ to $\Omega^0$ by $u^0 = u^{0, cutoff} + \sum_j \phi^1(\mathbf{r}^0)u^1_j$.
  \item We interpolate the potential from $\Omega^0$ to the off-grid atoms (with positions $\mathbf{r}_i$) by $u_i^{long} = \sum_j \phi^0(\mathbf{r}_i)u^0$.
  \item We calculate the potential for each atom $i$ by $U_i = 1/(4\pi\epsilon_0) (u_i^{short} + u_i^{long})$.
 \end{enumerate}
 
  Hardy precisely calculated the number of floating point operations needed for multilevel method \cite[ch.2]{Hardy2006}. It is
  \begin{equation}
  \begin{gathered}
  \label{eq:flop}
   \left( \frac{4}{3} \pi m + \frac{32}{3}\pi + \frac{81}{2} \right) \left( \frac{a}{h^*} \right)^3 N +\\ \left( 6p^3 + 31p^2 + 36p + 17 \right) N + \left( \left( \frac{4a}{h} \right)^3 + 14(p+2)\right) \frac{8}{7} \left( \frac{h^*}{h} \right)^3 N
   \end{gathered}
  \end{equation}
 where
 \begin{itemize}
  \item $h^* = N^{-1/3}L$ is the average distance between nearest neighbors 
  \item $h$ is the spacing in the finest grid 
  \item $a$  is the cutoff distance within which the interactions are considered short-range and calculated on current grid 
  \item $N$ is the number of atoms
  \item $L$ is the length of the box
  \item $\gamma$ is the smoothing function polynomial of degree $2m$
  \item $\phi$ is the basis function polynomial of degree $p$
 \end{itemize}
 
 In our simulation, $N=10^8$ and $L=500$\AA, therefore $h^* \approx 1$\AA. Usually \cite[in most evaluations]{Hardy2006}, $h=2$\AA~and the cutoff $a$ ranges from $8$\AA~to $12$\AA. When choosing $\gamma$ and $\phi$ functions, we need to consider the trade-off between the accuracy and the number of floating point operations. $C^2$ Taylor smoothing function ($m=2$) and cubic numerical Hermite interpolant ($p=3$) keep the accuracy in reasonable boundaries (less than 1\% relative error in potential for $>8$\AA~cutoff \cite[p.78]{Hardy2006}) while maintaining quite low cost. The number of floating point operations for these values is ($a=12$\AA)
 \begin{equation}
  77.7a^3N+566N+\frac{73a^3}{h^6}N+\frac{80}{h^3}N \approx 136813N \sim 10^{13}.
 \end{equation}

 \section{Parareal Method}
    Parareal method \cite{Lions2001, Gander2007, Maday2008} parallelizes the time domain by approximating the solution in time $t_1$ without accurate solution in time $t_0$ where $t_1 > t_0$. When solving differential equation, we seek the function of time $\mathbf{u}$ for initial condition $v$ in time $t>0$. 
%
  The exact solution is rarely known, more often we have precise enough approximation $\mathcal{F}_\tau(t, v)$ obtained by discretization with small timestep. We approximate $\{ \mathbf{u}(T_n) \}_n$ with
 \begin{equation}\{ \lambda_n = \mathcal{F}_{T_n-T_0}(T_0; v) = \mathcal{F}_{\Delta T_n}(T_{n-1}; \lambda_{n-1}) \}_n.\end{equation} where $T_i$ are evenly-spaced time points. We clearly see the sequential character of the problem, the solution in time $T_n$ can be calculated only after the solution in time $T_{n-1}$ is known.
 
 Parareal method proposes a sequence $\{ \lambda_n^k \}_n$ that converges to $\{ \lambda_n \}_n$ rapidly as $k \rightarrow \infty$ and that can be built in parallel. We introduce the second, coarse and less expensive approximation $\mathcal{G}_\tau(t, v)$. The sequence $\{ \lambda_n^k \}_n$ is then defined recursively as
 \begin{equation}\lambda_{n+1}^{k+1} = \mathcal{G}_{\Delta T}(T_n; \lambda_n^{k+1}) + \mathcal{F}_{\Delta T}(T_n; \lambda_n^k)-\mathcal{G}_{\Delta T}(T_n; \lambda_n^k)\end{equation}
 and it usually converges after a few iterations. The idea behind this sequence is to shift the inherent sequential nature of calculation from $\mathcal{F}$ to $\mathcal{G}$. We set initial, rough assessment of the result calculated only by $\mathcal{G}$ and then we can approximate the error -- by subtracting the results from precise calculation and from coarse calculation on the same rough inaccurate data -- and gradually improve the approximation of the result.
 
 Lets' take a more specific case. We have a molecular dynamics simulation with total length 1\ ps. We set the width of ``computational window''---within which we want to calculate the results parallel in time to e.g.~100\ fs. As the calculation proceeds (as explained in the next paragraph), we acquire converged results in first time points, we can shift the computational window to the right on time axis and continue to calculate results for later time points. 
 
\begin{figure}[h!]
  \caption{Computational flow of the parareal method.}
 \label{fig:parareal}
 \begin{center}
  \includegraphics[width=0.6\textwidth]{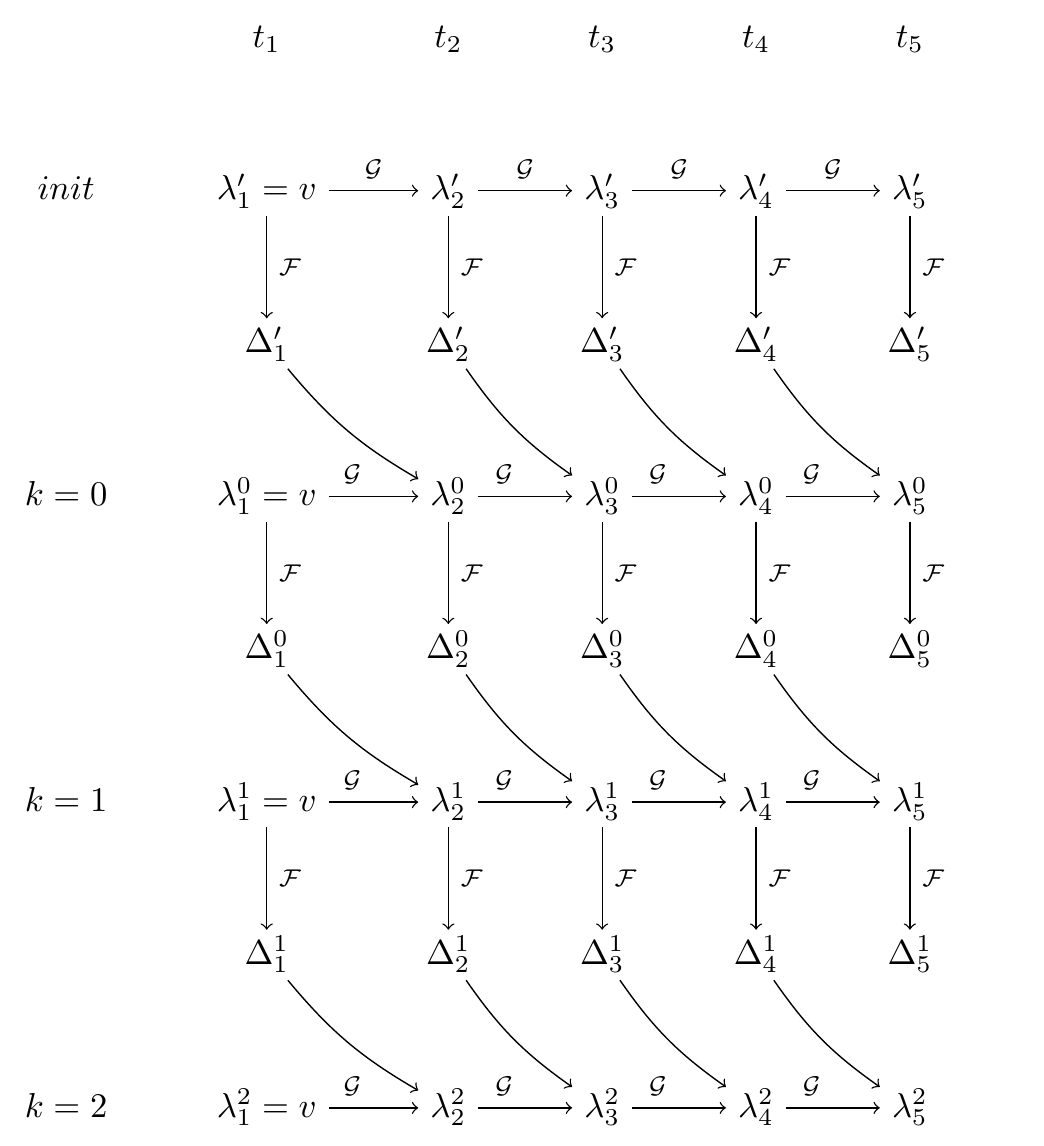}
  
\scriptsize
  $
\begin{array}{lllllll}
 & \lambda^2_1 = v &  & \lambda^2_3 \approx \lambda_3 \approx \mathcal{F}_{\Delta T}(T_2, \lambda_2) & & \lambda^2_5 \approx \lambda_5 \approx \mathcal{F}_{\Delta T}(T_4, \lambda_4) & \\
 \end{array}
$

$
   \begin{array}{lllllll}
  & &\lambda^2_2 \approx \lambda_2 \approx \mathcal{F}_{\Delta T}(T_1, v) & &  \lambda^2_4 \approx \lambda_4 \approx \mathcal{F}_{\Delta T}(T_3, \lambda_3) & &\\
 \end{array}
$
\normalsize
  \end{center}
  
 \end{figure}

 The figure \ref{fig:parareal} shows the parareal method that calculates $\lambda_n$ for five time points $n=1..5$ with the initial condition $v$. The results converge after three iterations $k=0,1,2$. We apply further notation for clarity:
 \begin{itemize}
  \item $f_n^k = \mathcal{F}_{\Delta T}(T_{n-1}, \lambda_{n-1}^k)$
  \item $g_n^k = \mathcal{G}_{\Delta T}(T_{n-1}, \lambda_{n-1}^k)$
  \item $\Delta_n^k = f_n^k-g_n^k$
  \item $'$ superscript replaces $(-1)^{st}$ iteration---initialization
 \end{itemize}
 
 The arrow shows the computational flow---what needs to be computed in what order and with what dependencies. In \textit{init} row, the calculation of $g'_n$ depends on already known $g'_{n-1}$. Downward arrows with $f_n^k$ and $\Delta_n^k$ represent parallel computation of precise approximation based on inaccurate data\footnote{Inaccurate data does not mean inaccurate initial condition, we emphasize that the \textit{coarse} approximation calculates results in time $t_n$ from results in $t_{n-1}$. This approximation gives us inaccurate data for calculation of result in the next step.}. In \textit{k} rows, $\lambda_n^k$ calculation depends on known $g_{n-1}^{k}$ and $\Delta_{n-1}^{k-1}$. After three iterations, we consider $\lambda_n^2$ very close to the result we would get by sequential computation of precise approximation based on accurate input data.
 
 \section{Multilevel Summation + Parareal Method}
  Our main research question focuses on how to combine multilevel method (that calculates the potential over spatial domain) and parareal time integration method (that approximates the potential a few timesteps ahead) into parareal multilevel summation method. We want to achieve computation parallel in both space and time. Multilevel summation method will calculate over spatial domain in parallel. The parareal method should ensure the parallelism in time. The classical MD method with multilevel approximation represents the fine approximation $\mathcal{F}$ of the parareal method. The main characteristic of the combination is the choice of parareal method's coarse approximation $\mathcal{G}$. Several approaches are available and listed below, the first four of them have obvious pitfalls therefore we want to develop the latter two.  For each of them, we compare them to classical molecular dynamics with MSM. Moreover, we examine two key aspects---the convergence and the ability to run in parallel. 
Considered approaches:
    \begin{itemize}
      \item discrete MD;
      \item longer timestep;
      \item coarse grid in MSM;
      \item shorter cutoff in MSM;
     \item simple cutoff method
     \item Wolf summation method
    \end{itemize}
    \paragraph{Discrete Molecular Dynamics} Coarse approximation $\mathcal{G}$ would run the simulation using another concept---instead of integrating the differential equations to capture the movements of atoms, we assume that they move with constant velocity unless a collision occurs \cite{Proctor2011}. Despite its very good convergence, discrete molecular dynamics would not work as it produces completely different trajectories than classical molecular dynamics. The overall properties of the system are similar, but the atoms move differently. As we are particularly interested how exactly the virus passes the membrane and how the atoms of virus and the membrane interact, this approach does not suit us.
    
    \paragraph{Longer Time Step} Coarse approximation $\mathcal{G}$ would have longer timestep in integration scheme, instead of $2\text{ fs}$ it would be e.g.~$10\text{ fs}$. The problem would appear in the beginning of the computation, even before $k=0$, when we need to calculate $\lambda'_{n+1} = \mathcal{G}(t_{n+1}, t_n, \lambda'_n)$ where $\lambda'_1 = v$. The coarse approximation $\mathcal{G}$ therefore needs to give at least a little reasonable results for every $t_n$ knowing only initial condition $v$ and its own result for $t_{n-1}$. Unfortunately, MD simulation with timestep larger or equal 5 fs give unusable results after just a few steps (in our own experiment with retinol and timestep 10 fs the system blew up after two steps). Half-converged results, e.g.~$\lambda_n^3$, instead of $\lambda'_n$, could save the convergence, but we would have to shorten the computational window. That would drastically reduce the ability to run in parallel.
    
   
    \paragraph{Coarse Grid} Coarse approximation $\mathcal{G}$ would run molecular dynamics with multilevel summation, but the finest grid would have spacing e.g.~four times larger than in the fine approximation $\mathcal{F}$. It would be easily implementable, however, the complexity assessment suggests it would not be efficient as the number of floating point operations depends mainly on the cutoff distance $a$, not on the finest grid spacing $h$.
    
    \paragraph{Shorter Cutoff Distance} Coarse approximation $\mathcal{G}$ would calculate the interactions on the current grid within shorter cutoff range than in the fine approximation $\mathcal{F}$ (e.g.~6\AA\ vs 12\AA\ for the finest grid)\footnote{Cutoff distance in MSM is the distance within which the part of potential is calculated on current grid.}. Bigger parts of the potential would be calculated on coarser grids, eventually large part of the potential would be calculated on the top level where it is cheaper than on finer grids. In evaluations \cite[ch.3-5]{Hardy2006}, Hardy examined the relative error with cutoff beginning at $8$\AA. Even with high order polynomials for smoothing functions and interpolants, the relative error in force ranges from 1\% to 5\%. Further analysis shows that the maximal theretical speedup is 7. As simple cutoff method or Wolf summation offer much larger theretical speedups, we abandon this possibility. 
    
    
     \paragraph{Simple Cutoff} Coarse approximation $\mathcal{G}$ would calculate long-range interactions by the simplest method available---sum up interactions with Coulomb equation between atoms within constant cutoff, usually set to 12~\AA. This method lacks the accuracy and even introduces some non-physical artefacts. However, the precise $\mathcal{F}$ function should correct them. The method is computationally cheap and rather accurate, therefore almost perfect candidate for $\mathcal{G}$ function.
    
    \paragraph{Wolf Summation Method} Coarse approximation would use Wolf summation method \cite{wolf99} to calculate electrostatics. This method stems from simple cutoff and with little additional computational cost it deals with most problems of simple cutoff method. Promising and worth trying.
        
    We suggest the Wolf summation method as far future work, for now lets' keep things simple and stay with simple cutoff method. To summarize, we want to run a simulation of molecular dynamics where we will parallelize in time the calculation of electrostatic potential by parareal method. Fine approximation $\mathcal{F}$ will calculate the potential by MSM with cutoff distance $a$, e.g.~12\AA. Coarse approximation $\mathcal{G}$ will calculate the potential by simple cutoff method. In both approximations, we will parallelize also the spatial domain. As one of the further improvements, we will use multiple-time-step integration method.
    
    \section{Current and Near Future Work}
    Now, we want to examine the theoretical speedup of the parareal method with MSM and simple cutoff along with the convergence of such method.
    
    There is another concept to explore. The parareal method offers an interesting possibility to even further shorten the time result. Many simulation of organic systems go though a metastable state during the simulation. The atoms do not move very much, the system waits for the peak in energy to cross a barrier or to get from local minimum. If we are able to assess (the order of) $\Delta^k_n$ very quickly without calculating $f_n^k$, we can detect that systems does not change much in the next step and electrostatics's evaluation done by $\mathcal{G}$ suffices. That way we can save the time and resources to calculate $f_n^k$ without much improvement in the result. Of course, the method itself, its accuracy and convergence conditions will need to be explored thoroughly.

\subsection{Complexity and ability to run in parallel}
Let's assume we are running a simulation with parareal multilevel summation method with $T_W$ time points in computational window, $T_{total} = W T_W$ time points in whole simulation and $K$ iterations of convergence. For the computation of $\lambda_n^{K-1}$ for all $n=1..T_W$ we need $T_W (K+1)-1$ calculations of function $\mathcal{G}$, $TK$ calculations of function $\mathcal{F}$ and $TK$ subtractions to get~$\Delta$. Assume we can calculate the function $\mathcal{G}$ with $Q_{\mathcal{G}}$ floating point operations on $P_{\mathcal{G}}$ processors in time $R_{\mathcal{G}}$ and the function $\mathcal{F}$ with $Q_{\mathcal{F}}$ floating point operations on $P_{\mathcal{F}}$ processors in time $R_{\mathcal{F}}$. If the number of processors depends only on spatial decomposition of calculation, then $P_{\mathcal{G}} \approx P_{\mathcal{F}}$ as we are simulating exactly the same system. For simplicity, lets' assume that 
the overhead of calculation caused by synchronization and awaiting communication should be similar both for $\mathcal{G}$ and $\mathcal{F}$. Therefore, we can assume 
\begin{equation}
\frac{R_{\mathcal{F}}}{R_{\mathcal{G}}} \approx \frac{Q_{\mathcal{F}}}{Q_{\mathcal{G}}}  = Q_{\mathcal{F}/\mathcal{G}}
\end{equation}

The cost of the calculation heavily depends on the distribution of tasks. In this context, a task represents the calculation of the function $\mathcal{G}$ or $\mathcal{F}$ on given input data $\lambda$ and $\Delta$. The distribution of tasks assigns the tasks to different processors. In preliminary analysis we came up with two distribution plans. Although they differ in number of used processors, the time speedup compared to the sequential algorithm is directly proportional to $Q_{\mathcal{F}/\mathcal{G}}$ in both.

\subsubsection{Calculation distribution plan 1}
The simplest plan assigns to computation unit P1\footnote{In these distribution plans, we consider a computation unit to have $P_{\mathcal{G}}$ or $P_{\mathcal{F}}$ processors.} calculations of $\mathcal{G}$ for all time points in computational window whereas $T_W = Q_{\mathcal{F}/\mathcal{G}}$. That should take $Q_{\mathcal{F}/\mathcal{G}}R_{\mathcal{G}} = R_{\mathcal{F}}$ time. After that, $Q_{\mathcal{F}/\mathcal{G}}$ units calculate in parallel $\mathcal{F}$ for all points in the computational window. After that, unit P1 calculates $\mathcal{G}$ for subsequent time points in the following computational window and continues analogously. We need only $Q_{\mathcal{F}/\mathcal{G}}$ units, but the time speedup is also very low.
\begin{equation}
\frac{T_{total}R_{\mathcal{F}}}{2\frac{T_{total}}{Q_{\mathcal{F}/\mathcal{G}}}R_{\mathcal{F}}} = \frac{Q_{\mathcal{F}/\mathcal{G}}}{2}
\end{equation}

\subsubsection{Calculation distribution plan 2}
\begin{figure}[tbh]
\caption{Calculation distribution plan 2.}
\label{fig:plan}
\begin{center}
\includegraphics[width=0.8\textwidth]{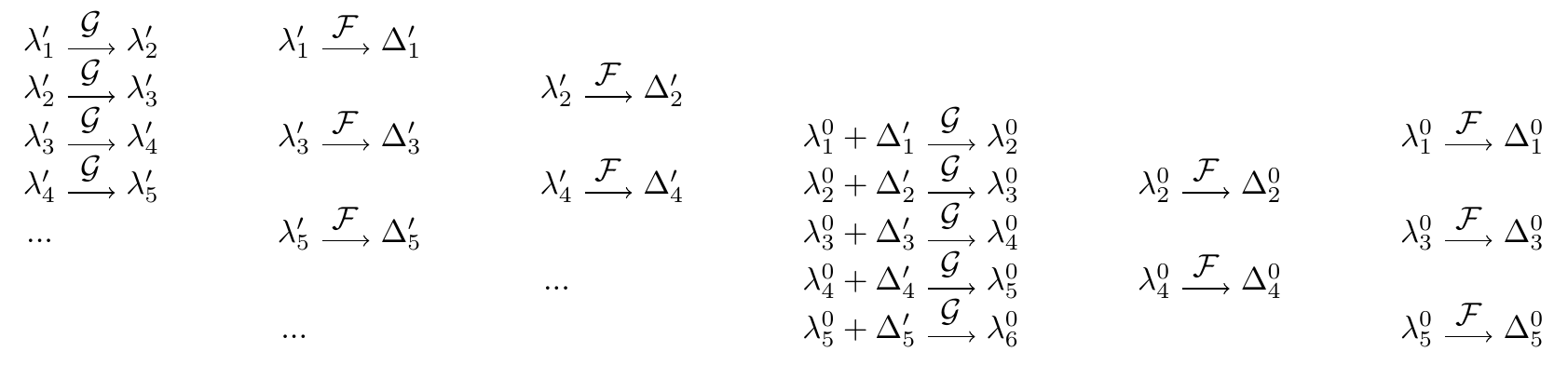}
\end{center}
\end{figure}

Figure \ref{fig:plan} depicts a more complex plan. It assigns the task of $\mathcal{F}$ calculation to the available (or new) unit as soon as it has the prerequisite - result of $\mathcal{G}$. We need $K$ units for $\mathcal{G}$ tasks and $KQ_{\mathcal{F}/\mathcal{G}}$ units for parallel calculation of $\mathcal{F}$. Time speedup is
\begin{equation}
 \frac{TR_{\mathcal{F}}}{(\frac{T}{Q_{\mathcal{F}/\mathcal{G}}}+K)R_{\mathcal{F}}} = \frac{Q_{\mathcal{F}/\mathcal{G}}}{1+\frac{K}{TQ_{\mathcal{F}/\mathcal{G}}}}
\end{equation}

\subsubsection{How to make it faster and more parallel?}
Now some specific numbers. Lets' assume $\mathcal{G}$ is simple cutoff method with cutoff 12\AA; $\mathcal{F}$ is MSM with $a_{\mathcal{F}}=12$\AA~and $h=2$\AA. Then
\begin{equation}
\begin{gathered}
Q_{\mathcal{G}} = 301N
Q_{\mathcal{F}} = 136813N\\
Q_{\mathcal{F}/\mathcal{G}} \approx 453\\
\end{gathered}
\end{equation}

Two calculation distribution plans presented above have time speedup proportional to $Q_{\mathcal{F}/\mathcal{G}}$. Perhaps some more advanced plan can increase that at the expense of increased number of processors needed. 

Now, we focus on finding the convergence conditions and implementation. 

\bibliographystyle{alpha}
\bibliography{library}

\newcommand{\etalchar}[1]{$^{#1}$}
\begin{thebibliography}{WKPE99}

\bibitem[BH86]{barneshut86}
J~Barnes and P~Hut.
\newblock {A hierarchical O(N log N) force-calculation algorithm}.
\newblock {\em Nature}, 324:446--449, 1986.

\bibitem[BHM00]{Briggs2000}
William~L Briggs, Van~Emden Henson, and Steve~F McCormick.
\newblock {\em {A Multigrid Tutorial}}.
\newblock SIAM Society for Industrial and Applied Mathematics, 2000.

\bibitem[DYP93]{darden93}
Tom Darden, Darrin York, and Lee Pedersen.
\newblock {Particle Mesh Ewald: An N.log(N) method for Ewald sums in large
  systems}.
\newblock {\em Journal of Chemical Physics}, 99:10089--10092, 1993.

\bibitem[EPa95]{essmann95}
U~Essmann, L~Perera, and Et~al.
\newblock {A smooth Particle Mesh Ewald Method}.
\newblock {\em Journal of Computational Physics}, 103(19):8577--8593, 1995.

\bibitem[Ewa21]{ewald21}
P~P Ewald.
\newblock {Die Berechnung optischer und elektrostatischer Gitterpotentiale}.
\newblock {\em Annalen der Physik}, 369(3):253--287, 1921.

\bibitem[GR87]{greengard87}
L~Greengard and V~Rokhlin.
\newblock {A Fast Algorithm for Particle Simulations}.
\newblock {\em Journal of Computational Physics}, 73:325--348, 1987.

\bibitem[GV07]{Gander2007}
Martin~J. Gander and Stefan Vandewalle.
\newblock {Analysis of the Parareal Time-Parallel Time-Integration Method}.
\newblock {\em SIAM Journal on Scientific Computing}, 29(2):556--578, January
  2007.

\bibitem[Har06]{Hardy2006}
David~Joseph Hardy.
\newblock {\em {Multilevel summation for the fast evaluation of forces for the
  simulation of biomolecules}}.
\newblock PhD thesis, 2006.

\bibitem[HSS09]{Hardy2009}
David~J Hardy, John~E Stone, and Klaus Schulten.
\newblock {Multilevel Summation of Electrostatic Potentials Using Graphics
  Processing Units.}
\newblock {\em Parallel computing}, 35(3):164--177, March 2009.

\bibitem[Jen07]{Jensen2007}
Frank Jensen.
\newblock {\em {Introduction to computational chemistry}}.
\newblock John Wiley \& Sons Ltd, Great Britain, 2 edition, 2007.

\bibitem[Koe06]{Koehl2006}
Patrice Koehl.
\newblock {Electrostatics calculations: latest methodological advances.}
\newblock {\em Current opinion in structural biology}, 16(2):142--51, April
  2006.

\bibitem[LMT01]{Lions2001}
Jacques-Louis Lions, Yvon Maday, and Gabriel Turinici.
\newblock {R\'{e}solution d'EDP par un sch\'{e}ma en temps «parar\'{e}el »}.
\newblock {\em Comptes Rendus de l'Acad\'{e}mie des Sciences - Series I -
  Mathematics}, 332(7):661--668, April 2001.

\bibitem[Mad08]{Maday2008}
Yvon Maday.
\newblock {The parareal in time algorithm}.
\newblock 2008.

\bibitem[PDD11]{Proctor2011}
Elizabeth Proctor, Feng Ding, and Nikolay~V. Dokholyan.
\newblock {Discrete molecular dynamics}.
\newblock {\em Wiley Interdisciplinary Reviews: Computational Molecular
  Science}, 1(1):80--92, January 2011.

\bibitem[SKE{\etalchar{+}}05]{shan05}
Yibing Shan, John~L Klepeis, Michael~P Eastwood, Ron~O Dror, and David~E Shaw.
\newblock {Gaussian split Ewald: A fast Ewald mesh method for molecular
  simulation}.
\newblock {\em Journal of Chemical Physics}, 122:54101, 2005.

\bibitem[ST82]{Stuben1982}
Klaus Stuben and Ulrich Trottenberg.
\newblock {Multigrid methods: Fundamental algorithms, model problem analysis
  and applications}.
\newblock In {\em Proceedings of the Conference Held at Koln-Porz}, volume 960,
  pages 1--176, 1982.

\bibitem[VMVF11]{Vacha2011}
Robert V\'{a}cha, Francisco~J. Martinez-Veracoechea, and Daan Frenkel.
\newblock {Receptor-Mediated Endocytosis of Nanoparticles of Various Shapes}.
\newblock {\em Nano Letters}, 11(12):5391--5395, 2011.

\bibitem[WKPE99]{wolf99}
D~Wolf, P~Keblinski, S~R Phillpot, and J~Eggebrecht.
\newblock {Exact method for the simulation of Coulombic systems by spherically
  truncated, pairwise r\^{}\{-1\} summation}.
\newblock {\em Journal of Chemical Physics}, 110(17):8254--8282, 1999.

\end{thebibliography}
      
\end{document}